\documentclass{amsart}%
\usepackage{amsmath}
\usepackage{amsfonts}
\usepackage{amssymb}
\usepackage{graphicx}%
\providecommand{\U}[1]{\protect\rule{.1in}{.1in}}
\newtheorem{theorem}{Theorem}[section]
\newtheorem{lemma}[theorem]{Lemma}
\theoremstyle{definition}
\newtheorem{definition}[theorem]{Definition}

\theoremstyle{remark}
\newtheorem{remark}[theorem]{Remark}
\numberwithin{equation}{section}

\begin{document}
\title[Real-Analyticity of Generalized Sine]{Real-Analyticity of Generalized Sine Functions with Two Parameters}
\author{Pisheng Ding}
\address{Department of Mathematics, Illinois State University, Normal, Illinois 61790}
\email{pding@ilstu.edu}
\date{}

\begin{abstract}
For any real number $p>1$ and any integer $n>1$, we identify the maximal real
interval on which the generalized sine function $\sin_{p,n}$ is real-analytic.
We achieve this by first proving that $\sin_{p,n}$ is analytic at $\frac{1}%
{2}\pi_{p,n}$ iff $p=m/(m-1)$ for some integer $m>1$, in which case we
determine the radius of convergence of the Taylor series at $\frac{1}{2}%
\pi_{p,n}$.

\end{abstract}
\maketitle

\section{Introduction}

For a real number $p>1$ and an integer $n>1$ such that $p$ and $n$ do not both
equal $2$, consider%
\[
F_{p,n}(y)=\int_{0}^{y}\frac{1}{(1-t^{n})^{1/p}}dt
\]
for $y\in\lbrack0,1]$. Denote $2F_{p,n}(1)$ by $\pi_{p,n}$. In the literature,
the generalized sine with parameters $(p,n)$, commonly denoted by $\sin_{p,n}%
$, is defined to be $\left(  F_{p,n}\right)  ^{-1}:\left[  0,\frac{1}{2}%
\pi_{p,n}\right]  \rightarrow\lbrack0,1]$. Clearly, $\sin_{p,n}$ is
real-analytic on $(0,\frac{1}{2}\pi_{p,n})$.

In this note, we find the maximal interval on which $\sin_{p,n}$ may be
analytically extended. Although we require that $p$ and $n$ do not both equal
2, the claims we make will turn out to be valid in the familiar classical case
$p=2=n$.

As $F_{p,n}$ is also monotonic on $(a_{n},1]$ where $a_{n}=-\infty$ for $n$
odd and $a_{n}=-1$ for $n$ even, $\sin_{p,n}$ has an obvious analytic
extension over the interval $(F_{p,n}(a_{n}),\frac{1}{2}\pi_{p,n})$. (However,
in much of the literature, $\sin_{p,n}$ is not extended in this manner for $n$
odd; see, e.g., \cite{Bh} and \cite{Ed}.) Henceforth, we consider the real
function%
\[
\sin_{p,n}:=F_{p,n}^{-1}:\left\{
\begin{tabular}
[c]{ll}%
$\left[  -\frac{1}{2}\pi_{p,n},\frac{1}{2}\pi_{p,n}\right]  \rightarrow
\lbrack-1,1]$, & for $n$ even\\
& \\
$\left(  F_{p,n}(-\infty),\frac{1}{2}\pi_{p,n}\right]  \rightarrow(-\infty
,1]$, & for $n$ odd
\end{tabular}
\ \ \ \right.  \text{.}%
\]
Note that $F_{p,n}$ and $\sin_{p,n}$ are odd iff $n$ is even.

To further extend $\sin_{p,n}$, a critical issue is the analyticity of
$\sin_{p,n}$ at $\frac{1}{2}\pi_{p,n}$. Concerning this, we shall prove the
following key result, in which we let $p^{\prime}:=p/(p-1)$, the so-called
conjugate of $p$, satisfying
\[
\frac{1}{p}+\frac{1}{p^{\prime}}=1\text{ ,}%
\]
and let%
\[
\mu(m,n):=\int_{1}^{\infty}\frac{1}{(t^{n}-1)^{(m-1)/m}}dt\text{ .}%
\]

\begin{theorem}
\label{Thm Analyticity at Pi/2}At $\frac{1}{2}\pi_{p,n}$, $\sin_{p,n}$ has an
analytic extension iff $p^{\prime}$ is an integer, in which case $\frac{1}%
{2}\pi_{p,n}$ is a critical point of order $p^{\prime}-1$ and the Taylor
series at $\frac{1}{2}\pi_{p,n}$ has radius of convergence $\mu(p^{\prime},n)$.
\end{theorem}

That is, $\sin_{p,n}$ is analytic at $\frac{1}{2}\pi_{p,n}$ iff $p=m/(m-1)$
for some integer $m>1$. With this proven, we will then focus on $\sin
_{m^{\prime},n}$ where $m>1$ is an integer.

In the proof of Theorem \ref{Thm Analyticity at Pi/2}, we apply
complex-analytic arguments to show that precisely when $p^{\prime}$ is an
integer does $\sin_{p,n}$ have a bounded analytic continuation on some
punctured disc $D(\frac{1}{2}\pi_{p,n};\epsilon)\setminus\{\frac{1}{2}%
\pi_{p,n}\}$. Furthermore, the proof also reveals how far the domain of
analyticity stretches in either direction along the real line; this is the
content of the next theorem, our main result, in which we let%
\[
\nu(m,n):=\int_{0}^{\infty}\frac{1}{(1+t^{n})^{(m-1)/m}}dt\text{ .}%
\]
Note that, for $n$ odd, $F_{p,n}(-\infty)=-\nu(p^{\prime},n)$.

\begin{theorem}
\label{Thm Maximal Interval}Suppose that $p^{\prime}$ is an integer.

\begin{enumerate}
\item If $p^{\prime}$ and $n$ are both even, then $\sin_{p,n}$ is analytic on
$%
\mathbb{R}
$.

\item If $p^{\prime}$ is even and $n$ is odd, then $\sin_{p,n}$ is analytic on
the interval%
\[
\left(  -\nu(p^{\prime},n),\;\pi_{p,n}+\nu(p^{\prime},n)\right)  \text{ .}%
\]

\item If $p^{\prime}$ is odd and $n$ is even, then $\sin_{p,n}$ is analytic on
the interval%
\[
\left(  -\frac{1}{2}\pi_{p,n}-\mu(p^{\prime},n),\;\frac{1}{2}\pi_{p,n}%
+\mu(p^{\prime},n)\right)  \text{ .}%
\]

\item If $p^{\prime}$ and $n$ are both odd, then $\sin_{p,n}$ is analytic on
the interval%
\[
\left(  -\nu(p^{\prime},n),\;\frac{1}{2}\pi_{p,n}+\mu(p^{\prime},n)\right)
\text{ .}%
\]

\end{enumerate}

Furthermore, each interval given is the maximal interval on which $\sin_{p,n}$
is real-analytic.
\end{theorem}

In fact, the endpoints of these maximal intervals in Theorem
\ref{Thm Maximal Interval} can all be expressed in terms of $\pi_{p,n}$, as
the following theorem claims.

\begin{theorem}
\label{Thm mu-nu-pi}Suppose that $p<n$ (which holds if $p^{\prime}$ is an
integer). Then,%
\[
\mu(p^{\prime},n)=\frac{\pi_{p,n}}{2\left(  \sin\frac{\pi}{p}\cot\frac{\pi}%
{n}-\cos\frac{\pi}{p}\right)  }\quad\text{and\quad}\nu(p^{\prime},n)=\frac
{\pi_{p,n}}{2\left(  \cos\frac{\pi}{n}-\sin\frac{\pi}{n}\cot\frac{\pi}%
{p}\right)  }\text{ .}%
\]

\end{theorem}

After reviewing in \S 2 some background matters concerning complex-analytic
continuation of $\sin_{p,n}$, we prove Theorem \ref{Thm Analyticity at Pi/2}
in \S 3, Theorem \ref{Thm Maximal Interval} in \S 4, and Theorem
\ref{Thm mu-nu-pi} in \S 5, in which we also apply our findings in a special case.

\section{Preliminaries}

It is shown in \cite{Ding} that $\sin_{p,n}$ has a complex-analytic
continuation on the interior of a closed polygon $\Pi_{p,n}\subset%
\mathbb{C}
$. We review certain details necessary for our development.

We first introduce notation. When we denote a point in $%
\mathbb{C}
$ by a capital letter, we write $AB$ for the (closed) line segment between $A$
and $B$, whereas we let $[A,B):=AB\setminus\{B\}$.

\begin{definition}
The following notations stay in effect throughout this note.

\begin{enumerate}
\item $A_{p,n}:=\frac{1}{2}\pi_{p,n}$; $J:=[1,\infty)$; $L_{n}:=\{te^{i\pi
/n}\mid t\geq0\}$; $W_{n}:=\{re^{i\theta}\mid r>0;\,\theta\in(0,\pi/n)\}$;
$W_{n}^{2}:=\{z^{2}\mid z\in W_{n}\}$.

\item For $p<n$, the closed triangle $\Delta_{p,n}$ has two of its three
vertices at $O:=0$ and $A_{p,n}$; its remaining vertex $C_{p,n}$ is such that%
\[
C_{p,n}\in L_{n}\quad\text{and}\quad\measuredangle OA_{p,n}C_{p,n}=\frac{\pi
}{p^{\prime}}\text{ .}%
\]

\item For $p\geq n$, $\Delta_{p,n}$ is the closed region in the first quadrant
bounded by the segment $OA_{p,n}$, the ray $L_{n}$, and the ray $\{\frac{1}%
{2}\pi_{p,n}+te^{i\pi/p}\mid t\geq0\}$.

\item With $\Delta_{p,n}^{\dag}$ denoting the image of $\Delta_{p,n}$ under
reflection across $L_{n}$,%
\[
\Pi_{p,n}:=\Delta_{p,n}\cup\Delta_{p,n}^{\dag}\text{ .}%
\]

\end{enumerate}
\end{definition}

We summarize some relevant facts that are immediate generalizations of the
main findings in \cite{Ding}.

\begin{lemma}
\label{Lem Key} (Cf. \cite{Ding}.)

\begin{enumerate}
\item $F_{p,n}:z\mapsto\int_{0}^{z}(1-\varsigma^{n})^{-1/p}d\varsigma$ is the
complex-analytic extension on $W_{n}^{2}$ of the real function $F_{p,n}$.

\item $F_{p,n}:W_{n}^{2}\rightarrow\mathring{\Pi}_{p,n}$ and $F_{p,n}%
:W_{n}\rightarrow\mathring{\Delta}_{p,n}$ are both conformal equivalences.

\item For $p\geq n$, the continuous extension of $F_{p,n}$ (also denoted by
$F_{p,n}$) is a homeomorphism $\partial W_{n}\rightarrow\partial\Delta_{p,n}$
with%
\[
F_{p,n}[L_{n}]=L_{n}\quad\text{and\quad}F_{p,n}\left[  J\right]  =\left\{
\frac{1}{2}\pi_{p,n}+te^{i\pi/p}\mid t\geq0\right\}  \text{ .}%
\]

\item For $p<n$, the continuous extension of $F_{p,n}$ (also denoted by
$F_{p,n}$) is a homeomorphism $\partial W_{n}\rightarrow\partial\Delta
_{p,n}\setminus\{C_{p,n}\}$ with%
\begin{align*}
F_{p,n}(1+t)  &  =\frac{1}{2}\pi_{p,n}+e^{i\pi/p}\int_{1}^{t}\frac{1}%
{(\tau^{n}-1)^{1/p}}d\tau\\
\text{and\quad}F_{p,n}(te^{i\pi/n})  &  =e^{i\pi/n}\int_{0}^{t}\frac
{1}{(1+\tau^{n})^{1/p}}d\tau
\end{align*}
for $t\geq0$. In particular,%
\[
F_{p,n}[L_{n}]=[O,C_{p,n})\quad\text{and}\quad F_{p,n}[J]=[A_{p,n}%
,C_{p,n})\text{ .}%
\]

\item The conformal equivalence $F_{p,n}^{-1}:\mathring{\Pi}_{p,n}\rightarrow
W_{n}^{2}$ analytically extends the real function $\sin_{p,n}$. (We also
denote this complex function $F_{p,n}^{-1}$ by $\sin_{p,n}$.)
\end{enumerate}
\end{lemma}

\section{Real-Analyticity of $\sin_{p,n}$ at $\pi_{p,n}/2$}

We now make the case for Theorem \ref{Thm Analyticity at Pi/2}.

Note that, by Lemma \ref{Lem Key}, the complex-analytic $\sin_{p,n}$ amplifies
angle at $A_{p,n}$ by the factor $p^{\prime}$.

If $\sin_{p,n}$ has a real-analytic extension on an open interval containing
$\frac{1}{2}\pi_{p,n}$, it will also have a complex-analytic extension on a
disc containing $A_{p,n}$, at which the angle-amplification factor $p^{\prime
}$ must be a positive integer, in which case $p=p^{\prime}/(p^{\prime}-1)<n$.
(Recall the assumption that $p$ and $n$ are not both $2$.)

On the other hand, suppose that $p^{\prime}$ is an integer.

Then $p=p^{\prime}/(p^{\prime}-1)<n$. Consider an open disc $D(A_{p,n};r)$
centered at $A_{p,n}$ of radius $r$. Choose $r$ so that $D(A_{p,n};r)\cap
\Pi_{p,n}$ is the circular sector $V_{p,n}(r)$ of radius $r$ with central
angle $\pi/p^{\prime}$ bounded by $A_{p,n}C_{p,n}$ and $A_{p,n}O$, i.e.,%
\[
D(A_{p,n};r)\cap\Pi_{p,n}=V_{p,n}(r):=\left\{  \frac{1}{2}\pi_{p,n}%
+te^{i\theta}\mid t\in\lbrack0,r);\,\theta\in\left[  \frac{\pi}{p},\pi\right]
\right\}  \text{.}%
\]
Since $\sin_{p,n}$ maps $A_{p,n}C_{p,n}$ and $A_{p,n}O$ into $%
\mathbb{R}
$, we can apply the Schwarz reflection principle repeatedly to analytically
extend $\sin_{p,n}$ (which is already analytic in $\mathring{V}_{p,n}(r)$) to
the $2p^{\prime}$ congruent sectors dividing $D(A_{p,n};r)$. This way,
$\sin_{p,n}$ is analytically continued on the punctured disk $D(A_{p,n}%
;r)\setminus\{A_{p,n}\}$. Since, at $A_{p,n}$, $\sin_{p,n}|_{\Delta_{p,n}}$ is
continuous, the Riemann principle of removable singularity implies that
$\sin_{p,n}$ is complex-analytic at $A_{p,n}$.

On the punctured disk $D(A_{p,n};r)\setminus\{A_{p,n}\}$, $\sin_{p,n}$ is
$p^{\prime}$-to-one, implying that $A_{p,n}$ is a critical point of order
$p^{\prime}-1$.

Let $R_{p,n}$ denote the radius of convergence of the Taylor series for
$\sin_{p,n}$ at $\frac{1}{2}\pi_{p,n}$. In light of Lemma \ref{Lem Key},
$|\sin_{p,n}z|\rightarrow\infty$ as $z\rightarrow C_{p,n}$ and consequently%
\[
R_{p,n}\leq|A_{p,n}C_{p,n}|=\mu(p^{\prime},n)\text{ .}%
\]
We now explain why this upper bound for $R_{p,n}$ is in fact its exact value.

When $p^{\prime}$ and $n$ both exceed $2$, it is a matter of plane geometry
that the circular sector $V_{p,n}(\mu(p^{\prime},n))\subset\Pi_{p,n}$.
Therefore, the analytic continuation of $\sin_{n,p}$ constructed above is
defined on $D(A_{p,n};\mu(p^{\prime},n))$. Hence, in this case, $R_{p,n}%
=\mu(p^{\prime},n)$.

If $(p^{\prime},n)=(2,3)$ or $(p^{\prime},n)=(3,2)$, then $V_{p,n}%
(\mu(p^{\prime},n))\nsubseteq\Pi_{p,n}$. In both cases, $\Delta_{p,n}$ is a
$30^{\circ}$-$60^{\circ}$-$90^{\circ}$ triangle; $\Pi_{2,3}$ is the union of
two copies of such a triangle joined along the hypotenuse, whereas
$\Pi_{3^{\prime},2}$ is an equilateral triangle. In each case, the mapping
properties of $\sin_{p,n}$ (as outlined in Lemma \ref{Lem Key}) allows
$\sin_{p,n}$ to be further analytically continued (by the Schwarz reflection
principle) on $V_{p,n}(\mu(p^{\prime},n))\setminus\Pi_{p,n}$; key to this
argument is the fact that every point in $V_{p,n}(\mu(p^{\prime}%
,n))\setminus\Pi_{p,n}$ has an image in $\Pi_{p,n}$ under reflection across a
side of $\Pi_{p,n}$. Thus, $\sin_{p,n}$ is again analytic on $D(A_{p,n}%
;\mu(p^{\prime},n))$, as in the previous case, and $R_{p,n}=\mu(p^{\prime},n)$.

Finally, note that the formula for $R_{p,n}$ remains valid for $p=2=n$, as
$R_{2,2}=\infty=\mu(2,2)$.

\section{Maximal Interval of Real-Analyticity}

The analytic continuation of $\sin_{p,n}$ on the disc $D(A_{p,n};\mu
(p^{\prime},n))$ constructed in \S 3 provides explicit description of
$\sin_{p,n}x$ for $x>\frac{1}{2}\pi_{p,n}$. To explicate this, we first make
an observation.

\begin{lemma}
\label{Lem Rotational Invariance}Suppose that $p^{\prime}$ is an integer. Let
$\rho$ be rotation around $A_{p,n}$ by angle $2\pi/p^{\prime}$. Then, for
$z\in D(A_{p,n};\mu(p^{\prime},n))$,%
\[
\sin_{p,n}\rho(z)=\sin_{p,n}z\text{ .}%
\]

\end{lemma}

\begin{proof}
For each $k\in%
\mathbb{Z}
$, let $\ell_{k}:=\left\{  \frac{1}{2}\pi_{p,n}+te^{ik\pi/p^{\prime}}\mid
t\in\lbrack0,\mu(p^{\prime},n))\right\}  $, which is a radius of
$D(A_{p,n};\mu(p^{\prime},n))$ and is mapped by $\sin_{p,n}$ into $%
\mathbb{R}
$ (as shown in the preceding section). These $2p^{\prime}$ radii divide
$D(A_{p,n};\mu(p^{\prime},n))$ into $2p^{\prime}$ (relatively closed) sectors.
Let $\sigma_{k}$ denote reflection across $\ell_{k}$. If $z\in D(A_{p,n}%
;\mu(p^{\prime},n))$ is in the sector bounded between $\ell_{k}$ and
$\ell_{k+1}$, then, by the Schwarz reflection principle,%
\[
\sin_{p,n}\left[  \left(  \sigma_{k+1}\circ\sigma_{k}\right)  (z)\right]
=\sin_{p,n}z\text{ .}%
\]
On the other hand, $\sigma_{k+1}\circ\sigma_{k}=\rho$. The claimed identity
then follows.
\end{proof}

We now argue Theorem \ref{Thm Maximal Interval}. There are two cases
corresponding to the parity of $p^{\prime}$; each case in turn has two
subcases corresponding to the parity of $n$.

\begin{description}
\item[Case 1] $p^{\prime}$ is even. For $t\in(0,\mu(p^{\prime},n))$, by Lemma
\ref{Lem Rotational Invariance},%
\begin{equation}
\sin_{p,n}\left(  \frac{1}{2}\pi_{p,n}+t\right)  =\sin_{p,n}\left(
\rho^{p^{\prime}/2}\left(  \frac{1}{2}\pi_{p,n}+t\right)  \right)  =\sin
_{p,n}\left(  \frac{1}{2}\pi_{p,n}-t\right)  \text{.} \label{p' even}%
\end{equation}
This identity allows interpretation of $\sin_{p,n}\left(  \frac{1}{2}\pi
_{p,n}+t\right)  $ for those positive $t$ for which $\sin_{p,n}\left(
\frac{1}{2}\pi_{p,n}-t\right)  $ is meaningful. We discuss two subcases.

\begin{description}
\item[Subcase 1] $n$ is odd. In (\ref{p' even}), we may allow all positive $t$
with%
\[
t<\frac{1}{2}\pi_{p,n}-F_{p,n}(-\infty)=\frac{1}{2}\pi_{p,n}-\nu(p^{\prime
},n)
\]
and thus $\sin_{p,n}$ is defined and analytic on $\left(  -\nu(p^{\prime
},n),\;\pi_{p,n}-\nu(p^{\prime},n)\right)  $. This interval is maximal, as
$\sin_{n,p}x\rightarrow-\infty$ at the two endpoints.

\item[Subcase 2] $n$ is even. Appealing to the oddness of $\sin_{p,n}$ in
conjunction with (\ref{p' even}), we may extend $\sin_{p,n}$ around $-\frac
{1}{2}\pi_{p,n}$ by requiring that%
\begin{equation}
\sin_{p,n}\left(  -\frac{1}{2}\pi_{p,n}-t\right)  =\sin_{p,n}\left(  -\frac
{1}{2}\pi_{p,n}+t\right)  \text{.} \label{p' even - n even}%
\end{equation}
On $%
\mathbb{R}
$, the composition of reflections about $-\frac{1}{2}\pi_{p,n}$ and $\frac
{1}{2}\pi_{p,n}$ is translation by $2\pi_{p,n}$. Thus, the conjunction of
(\ref{p' even}) and (\ref{p' even - n even}) implies that $\sin_{p,n}$ can be
periodically extended over $%
\mathbb{R}
$:%
\[
\sin_{p,n}(x+2\pi_{p,n})=\sin_{p,n}x\text{ .}%
\]
This periodic extension is clearly analytic on $%
\mathbb{R}
$.
\end{description}

\item[Case 2] $p^{\prime}$ is odd. For $t\in(0,\mu(p^{\prime},n))$,%
\[
\sin_{p,n}\left(  \frac{1}{2}\pi_{p,n}+t\right)  =\sin_{p,n}\left(
\rho^{(p^{\prime}-1)/2}\left(  \frac{1}{2}\pi_{p,n}+t\right)  \right)
\text{.}%
\]
Since $\rho^{(p^{\prime}-1)/2}$ is rotation by angle $(1-(1/p^{\prime}%
))\pi=\pi/p$,%
\[
\rho^{(p^{\prime}-1)/2}\left(  \frac{1}{2}\pi_{p,n}+t\right)  =\frac{1}{2}%
\pi_{p,n}+te^{i\pi/p}%
\]
and therefore%
\begin{equation}
\sin_{p,n}\left(  \frac{1}{2}\pi_{p,n}+t\right)  =\sin_{p,n}\left(  \frac
{1}{2}\pi_{p,n}+te^{i\pi/p}\right)  \label{p' odd}%
\end{equation}
for $t\in(0,\mu(p^{\prime},n))$. We again address two subcases.

\begin{description}
\item[Subcase 1] $n$ is odd. With (\ref{p' odd}), $\sin_{p,n}$ is defined and
analytic on $\left(  -\nu(p^{\prime},n),\;\frac{1}{2}\pi_{p,n}+\mu(p^{\prime
},n)\right)  $. As in Subcase 1 of Case 1, this interval is obviously maximal.

\item[Subcase 2] $n$ is even. With (\ref{p' odd}), $\sin_{p,n}$ is already
defined on%
\[
\left(  -\frac{1}{2}\pi_{p,n},\;\frac{1}{2}\pi_{p,n}+\mu(p^{\prime},n)\right)
\text{ .}%
\]
Appealing to the oddness of $\sin_{p,n}$, we extend $\sin_{p,n}$ to%
\[
\left(  -\frac{1}{2}\pi_{p,n}-\mu(p^{\prime},n),\;\frac{1}{2}\pi_{p,n}%
+\mu(p^{\prime},n)\right)  \text{ .}%
\]
As $|\sin_{p,n}x|\rightarrow\infty$\ as $x$ approaches either endpoint, this
interval is maximal.
\end{description}
\end{description}

\begin{remark}
Suppose that $p^{\prime}$ is odd. Then we may describe $\sin_{p,n}(\frac{1}%
{2}\pi_{p,n}+t)$ in \textit{real} (as opposed to \textit{complex}) terms. Let%
\[
G_{p,n}(y)=\int_{1}^{y}\frac{1}{(t^{n}-1)^{1/p}}dt\text{ .}%
\]
Clearly, $G_{p,n}$ is invertible on $(1,\infty)$. We then have%
\[
\sin_{p,n}(\frac{1}{2}\pi_{p,n}+t)=G^{-1}(t)\text{ .}%
\]

\end{remark}

\section{Relations Among $\pi_{p,n}$, $\mu(p^{\prime},n)$, and $\nu(p^{\prime
},n)$}

Considering $\Delta_{p,n}$ where $p<n$, we note the following relation among
the three numbers $\mu(p^{\prime},n)$, $\nu(p^{\prime},n)$, and $\pi_{p,n}$.

\begin{lemma}
Suppose that $p<n$. The numbers $\mu(p^{\prime},n)$, $\nu(p^{\prime},n)$,
$\frac{1}{2}\pi_{p,n}$ are the side lengths of a triangle whose interior
angles opposite the three sides are of measure $\pi/n$, $\pi/p^{\prime}$,
$(\frac{1}{p}-\frac{1}{n})\pi$, respectively. As a consequence,%
\[
\nu(p^{\prime},n)\sin\frac{\pi}{n}=\mu(p^{\prime},n)\sin\frac{\pi}{p^{\prime}%
}\text{ ,}%
\]
and%
\[
\nu(p^{\prime},n)\cos\frac{\pi}{n}+\mu(p^{\prime},n)\cos\frac{\pi}{p^{\prime}%
}=\frac{1}{2}\pi_{p,n}\text{ .}%
\]

\end{lemma}

Theorem \ref{Thm mu-nu-pi} follows at once upon noting that $\cos\frac{\pi
}{p^{\prime}}=-\cos\frac{\pi}{p}$ and $\sin\frac{\pi}{p^{\prime}}=\sin
\frac{\pi}{p}$.

\begin{remark}
By Theorem \ref{Thm mu-nu-pi},%
\[
\mu(n,n)=\frac{1}{4}\pi_{n^{\prime},n}\sec\frac{\pi}{n}=\nu(n,n)\text{ .}%
\]
When $p^{\prime}=n$, Theorem \ref{Thm Maximal Interval} reduces to two cases:

\begin{enumerate}
\item If $n$ is even, $\sin_{n^{\prime},n}$ is analytic and periodic on $%
\mathbb{R}
$. This reflects the fact that $\sin_{n^{\prime},n}t$ is the $y$-coordinate of
a point on the so-called $n$-circle $\{(x,y)\mid x^{n}+y^{n}=1\}$ where $t$ is
twice the area of the corresponding \textquotedblleft
circular\textquotedblright\ sector.

\item If $n$ is odd, the maximal interval on which $\sin_{n^{\prime},n}$ is
analytic is%
\[
\left(  -\frac{1}{4}\pi_{n^{\prime},n}\sec\frac{\pi}{n},\;\frac{1}{2}%
\pi_{n^{\prime},n}+\frac{1}{4}\pi_{n^{\prime},n}\sec\frac{\pi}{n}\right)
\text{.}%
\]
The length of this interval, i.e., $\frac{1}{2}\pi_{n^{\prime},n}(1+\sec
\frac{\pi}{n})$, is twice the area bounded between the curve $\{(x,y)\mid
x^{n}+y^{n}=1\}$ and its asymptote. When $n=3$, this interval is $\left(
-\frac{1}{2}\pi_{\frac{3}{2},3},\,\pi_{\frac{3}{2},3}\right)  $, whose length
$\frac{3}{2}\pi_{\frac{3}{2},3}$ is a fundamental period of the elliptic
function $\sin_{\frac{3}{2},3}$.
\end{enumerate}

For an inviting exposition on this geometric interpretation of $\sin
_{n^{\prime},n}$, see \cite{D-C}.
\end{remark}

\end{document}